%% file: main.tex
\documentclass[11pt]{amsart}

\input{packages}
\input{theorems}

\input{definitions}

\usepackage{bmpsize}

\begin{document}

\title{Planar reinforced $k$-out percolation}

\author{Gideon Amir}
\address[Gideon Amir]{Department of Mathematics,  Bar-Ilan University, Ramat-Gan, 5290002 Israel}
\email{gideon.amir@biu.ac.il}
\author{Markus Heydenreich}
\address[Markus Heydenreich]{Institut für Mathematik, Universität Augsburg, 86135 Augsburg, Germany}
\email{markus.heydenreich@uni-a.de}
\author{Christian Hirsch}
\address[Christian Hirsch]{Department of Mathematics, Aarhus University, Ny Munkegade 118, 8000 Aarhus C, Denmark}
\address[Christian Hirsch]{DIGIT Center, Aarhus University, Finlandsgade 22, 8200 Aarhus N, Denmark}
\email{hirsch@math.au.dk}

\begin{abstract}
 We investigate the percolation properties of a planar reinforced network model. In this model, at every time step, every vertex chooses $k \ge 1$ incident edges, whose weight is then increased by 1. The choice of this $k$-tuple occurs proportionally to the product of the corresponding edge weights raised to some power $\a > 0$.

	Our investigations are guided by the conjecture that the set of infinitely reinforced edges percolates for $k = 2$ and $\a \gg 1$. First, we study the case $\a = \ff$, where we show the percolation for $k = 2$ after adding arbitrarily sparse independent sprinkling and also allowing dual connectivities. We also derive a finite-size criterion for  percolation without sprinkling. Then, we extend this finite-size criterion to the $\a < \ff$ case. Finally, {we verify these conditions numerically.}
\end{abstract}
\maketitle

\input{intro}
\input{model}

\input{ind}

\input{osss}

\input{reinf}
\input{sim}

\section*{Acknowledgments}
This work was initiated during the workshop \emph{Adaptive Learning and Opinion Dynamics in Social Networks} at Bar-Ilan University in Tel Aviv.
%

%%%%%%%%%%%%%%%%%%%%%%%%%
\bibliographystyle{abbrv}
\bibliography{lit}
\end{document}

%% file: packages.tex
\usepackage[dvipsnames]{xcolor}
\usepackage{tikz}
\usepackage{pgf}
\usetikzlibrary{calc}
\usetikzlibrary{patterns}
\usetikzlibrary{arrows}
\usetikzlibrary{decorations.pathreplacing}

\usepackage{amsmath,bm,bbm,amsthm, amssymb}
\usepackage{mathtools}
\usepackage{fullpage}
\usepackage{color}
\usepackage{dsfont}
\usepackage{animate}
\usepackage{etoolbox}
\usepackage{comment}
\usepackage[utf8]{inputenc}
\usepackage{pgfplots}

\usepackage{tcolorbox}
\usepackage{adjustbox}
\usepackage[final]{pdfpages}
\usepackage[ruled,vlined,linesnumbered,algosection,resetcount]{algorithm2e}
\usepackage{blkarray}
\usepackage{arydshln}
\usepackage{wrapfig}
\usepackage[colorlinks=true]{hyperref}
\usepackage{todonotes}

%% file: theorems.tex
\theoremstyle{plain}
\newtheorem{theorem}{Theorem}%[chapter]
\newtheorem{proposition}[theorem]{Proposition}
\newtheorem{corollary}[theorem]{Corollary}
\newtheorem{lemma}[theorem]{Lemma}

\theoremstyle{definition}

\newtheorem{conjecture}[theorem]{Conjecture}

\theoremstyle{remark}

%% file: definitions.tex
\def\th{\theta}
\def\Z{\mathbb Z}

\def\P{\mathbb P}

\def\E{\mathbb E}

\def\00{\mathbf 0}

\def\bet{\begin{theorem}}
\def\ent{\end{theorem}}
\def\bec{\begin{corollary}}
\def\enc{\end{corollary}}
\def\bep{\begin{proof}}
\def\enp{\end{proof}}
\def\f{\frac}
\def\om{\omega}
\def\a{\alpha}

\def\es{\emptyset}
\def\su{\subseteq}
\def\ms{\mathsf}
\def\co{\colon}
\def\N{\mathbb N}
\def\mc{ \mathcal}
\def\ff{\infty}

\def\EE{\mc E}

\def\one{\mathds1}
\def\de{\delta}
\def\d{\mathrm d}

\renewcommand\geq{\geqslant}
\renewcommand\le{\leqslant}
\renewcommand\ge{\geqslant}

\def\e{\varepsilon}
\def\sp{{\ms{spr}}}
\def\ps{p_{\sp}}
\def\Inf{\ms{Inf}}
\def\lrsa{\leftrightsquigarrow}
\def\bel{\begin{lemma}}
\def\enl{\end{lemma}}
\def\im{\item}
\def\been{\begin{enumerate}}
\def\enen{\end{enumerate}}
\def\beit{\begin{itemize}}
\def\enit{\end{itemize}}
\def\pa{\partial}
\def\sm{\setminus}

\def\k_d{\kappa_d}

\def\tff{\uparrow\infty}
\def\s{\sigma}
\def\co{\colon}
\def\ti{\times}
\def\k{\kappa}

\def\bepr{\begin{proposition}}
\def\enpr{\end{proposition}}

\def\ba{\,|\,}

\def\Var{\ms{Var}}

\def\eff{E_{\ms{fs}, \ff}}
\def\efa{E_{\ms{fs}, \a}}
\def\efap{E_{\ms{fs}, \a, p_* }}

%% file: intro.tex
\section{Introduction}
\label{sec:int}

%
%REINF
%
Mathematical processes with reinforcements are a vibrant research field as they appear in a wide range of applications. For instance, in economics, dominant companies that already have a large market share are more likely to expand their share even further. In neuroscience, synapses that have been used successfully in the past are more likely to be strengthened in the future. 
This effect, known as neuroplasticity, is the reason for a dynamical evolution of the neural network. An extreme form has been suggested by Markram et al.~as \emph{tabula-rasa hypothesis} \cite{markram}. 

The challenge for mathematicians is to explain how local reinforcement mechanisms drive global network properties. This has initiated the novel research field of graph-based reinforcement models. Here, we address the question whether it is  possible to design a reinforcement mechanism that gives rise to a network structure that is sparse but can still connect distant areas in the brain.

%
%MODEL DEF
%
In this work, we propose a percolating model of graph-based reinforcement, which is based on the idea of \emph{$k$-out networks}. {The evolution of the network takes place in rounds.} 
In each round, each vertex $z\in \Z^d$ reinforces precisely $k$ of its incident edges. The selection of such a $k$-tuple is proportional to the product of the current edge weights raised to some reinforcement parameter $\a > 0$. A particularly important special setting is the case where $\a = \ff$, where even a small weight difference makes an edge infinitely stronger than another one with lower weight. 

%
%MAIN THEORY RESULT
%
The question of percolation of the infinitely-reinforced edges is the main focus of our interest. More specifically, we concentrate on the case $d = k = 2$. Here, our overarching conjecture is that percolation occurs once $\a$ is sufficiently large. Our main result, Theorem \ref{thm:osss} below, is a major step towards this conjecture in the case $\a = \ff$. More precisely, we show that there is percolation after an arbitrarily sparse sprinkling of independent open edges when using a modified notion of connectivity that will be specified precisely below. Moreover, we provide several rigorously proven finite-range criteria for percolation, for which we then obtain overwhelming numerical evidence by Monte Carlo simulation. All of these results are stated in full detail in Section \ref{sec:mod} below.

%
%RELATION
%
We conclude the present introduction by elaborating on the connections of our work to two vibrant fields of research, namely i) graph-based reinforcement, and  ii) dependent models for percolation, especially those violating the FKG condition.
\bigskip

%
%GRAPH
%
\noindent {\bf Graph-based reinforcement.}
Pioneering works in the area of graph-based reinforcment are \cite{benaim,lucas} where the reinforcement is set on the vertices of the graph. In the context of neuroscience, it is more natural to consider edge-based reinforcement, and those were introduced in the seminal work \cite{reinf1}. One of the major drawback of this model is that even for $\a \gg 1$, it does not give rise to percolating structures \cite{reinf2}. Percolation could be achieved when relying on tree-based models \cite{ext1,ext2}. However, this imposes already very strong a priori structures. An interesting alternative form of graph reinforcement comes from ant-based walks \cite{ant1,ant2}. 
\bigskip

\noindent{\bf Percolation without FKG.}
In recent years, the investigation of dependent percolation models has attracted considerable attention. We refer the reader to \cite{voronoi} and references therein. However, although a broad range of models has been considered, the vast majority satisfies the \emph{FKG inequality}, i.e., positive correlation of increasing events. When considering questions beyond the validity of the FKG inequality, only a few isolated results are available. The results treated in \cite{muir,muir2} fundamentally rely on the Gaussian field and Poisson point process setting, respectively, and therefore do not apply in our setting. To solve this problem in our setting, we develop a completely novel stochastic domination property.

While our investigation concern a reinforced model for $k$-out percolation, we stress here that already the question of \emph{independent} $k$-out percolation has been the topic of vigorous research. More precisely, to put our investigation into perspective, we first review the independent $k$-out model discussed in \cite{koppl}. Here,  at each site $z \in \Z^d$, we select independently and uniformly at random a $k$-element subset $E_z$ of all bonds incident at $z$. We let $G_k^d = \bigcup_{z \in \Z^d}E_z$ denote the random graph obtained as the union of all these bonds. Then, we let
$$k_c(d) := \min\{k \ge 1\co \P(\text{$G_k^d$ percolates}) > 0\}$$
denote the smallest value of $k \ge 1$ such that $G_k^d$ percolates, noting that $k_c(d) \le d + 1$. A key property of this model is the negative correlation of the vacant edges. 
%
%THM IND
%
\bepr[Percolation of independent $d$-out percolation on $\Z^d$; \cite{koppl}]
\label{thm:ind}
For every $d \ge 3$, we have that $k_c(d) \le 3$. Moreover, $k_c(2) = 2$. 
\enpr

The rest of the manuscript is organized as follows. In Section \ref{sec:mod}, we give a detailed description of our $k$-out percolation model and state all our main results. 
Section \ref{sec:dom} contains the proof of a crucial domination property. 
In Section \ref{sec:osss}, we prove a sprinkling result by establishing sharp phase transition and combining it with the domination result. 
Subsequently, in Section \ref{sec:reinf}, we derive certain finite-size criteria for percolation, for which we then present overwhelming numerical evidence in Section \ref{sec:sim}.

%% file: model.tex
%
%SEC MOD
%
\section{Model and main results}
\label{sec:mod}

In this work, we consider a reinforced model for $k$-out percolation, which depends on a reinforcement parameter $\a > 0$. First, in Section \ref{ssec:ff}, we focus on the case $\a = \ff$. Then, in Section \ref{ssec:fin}, we consider the setting of finite $\a$.

%
%A FF
%
\subsection{The case $\a = \ff$}
\label{ssec:ff}

In this model, we assign weights $(W_e)$ to the edges $E$ of  the hyper-cubic lattice $\Z^d = (V, E)$, which evolve in discrete time steps. Initially, all edges have weight $W_e(0) = 1$. Then, at every discrete time $t \ge 1$ the weights are updated as follows.  
\been
\im 
At every site $z \in V$, we select one tuple $\s \in E_z^{(k)}$, where $E_z^{(k)}$ denote the family of $k$-tuples of edges incident to $z$. More precisely, we order all adjacent edges $e_1,\dots,e_{2d}$ according to their momentary weights $W_{e_1}(t),\dots,W_{e_{2d}}(t)$ in decreasing order (breaking ties by independent coin flips). We set $\s$ to be the first $k$ of these edges in this decreasing order.
\im  The weight $W_e$ of each edge $e\in E$ is then increased by the number of times it was selected. This selection is done simultaneously for all vertices $z\in V$ based on the weights $\{W_e(t)\}_{e\in E}$. As it can happen that an edge is selected by each of its two endpoints, one has $W_e(t+1)-W_e(t)\in\{0,1,2\}$.

\enen

This results in the ``final'' configuration
\begin{align}
	\label{eq:eeff}
\EE_\ff := \big\{ e\in E\co \liminf_{t\tff}(W_e(t) -W_e(t - 1))\ge 1 \big\} 
\end{align}
of  edges that are eventually reinforced in every round. Our arguments in Section \ref{sec:dom} show that this set coincides with the set of edges that are reinforced infinitely often. The same set of arguments shows that our percolation  model is stochastically dominated by the  independent $k$-out model mentioned in Proposition \ref{thm:ind}. 
This is explained in the proof of Proposition \ref{pr:dom} below.

To simplify the terminology, we call an edge $e \in E$ \emph{open} if $e \in \EE_\ff$ and \emph{closed} otherwise.
As elucidated in the introduction, we hypothesize that $\mc E_\ff$ percolates for $d = k = 2$. That is, we make the following conjecture.

\begin{conjecture}[Percolation for $d = k = 2$]
\label{conj:reinf}
Assume that $d= k = 2$.  Then,
	$\P(\EE_{\ff} \text{ percolates}) = 1.$
\end{conjecture}

The aim of this work is to present a series of results
supporting this core conjecture. First,  superpose the edges of interest it with some \emph{sprinkling}, i.e., by adding an independent Bernoulli percolation process with some small parameter $\e > 0$. We write $\EE_{\e, \ms{spr}}$ for the sprinkled set of bonds and put $\EE_{\ff, \e} := \EE_\ff \cup \EE_{\e, \ms{spr}}$. Adding the sprinkled bonds makes it substantially easier to find an infinite connected component. 
One of the striking properties of planar Bernoulli bond percolation is its self-duality. While our dependent model is still planar, it is no longer self-dual. This means that when applying classical techniques from percolation theory, the implications are often different than in Bernoulli bond percolation.  However, for our argument to work, we still need to allow the alternative that there is an infinite component where the adjacency notation from the original lattice is replaced by the adjacencies from the dual lattice. That is, two edges of $\Z^2$ are dually adjacent if there exists $z \in \Z^2$ such that $\partial (z +[0, 1]^2)$ contains both edges.

{We say that the model percolates if there exists an infinite sequence of distinct edges $(e_0,e_1,e_2,\dots)$ in $\EE_{\ff}$ such that for all $i\in\N$, we have that $e_{i-1}$ and $e_i$ are adjacent. 
We further say that the model \emph{percolates dually} if there exists an infinite sequence of distinct edges $(e_0,e_1,e_2,\dots)$ in $\EE_{\ff}$ such that for  all $i\in\N$, we have that $e_{i-1}$ and $e_i$ are \emph{dually adjacent}. }

%
%SPR
%
\bet[Percolation of sprinkled  percolation]
\label{thm:osss}
 Let $d= k = 2$ and $\e > 0$. Then, 
 $$\P(\EE_{\ff, \e}\text{ percolates}) = 1\; \text{ or }\; \P(\EE_{\ff, \e}\text{ percolates dually}) = 1.$$ 
\ent

%	
%PRF OVERVIEW
%
The proof of Theorem \ref{thm:osss} relies fundamentally on the OSSS inequality from \cite{osss}. In order to apply it, we need to derive bounds on influences and revealment probabilities. Here, we note that the influences are intimately related to probabilities of pivotal events, which have been studied for the geometrically challenging models such as Voronoi  percolation \cite{dembin}. However, the arguments of \cite{dembin} heavily rely on the FKG inequality, which is not available for our model. While recently, there has been progress in applying the OSSS inequality in models without FKG inequality \cite{muir,muir2}, these arguments are highly model dependent and do not extend to the present setting. We will address this challenge through a completely novel stochastic domination property in our setting, see Proposition \ref{pr:dom} below. It is when extending the techniques of \cite{muir,muir2} to our model that the finite range of dependance will be crucial.

Conjecture \ref{conj:reinf} is a percolation result for a model with parameter $k$, which is similar to the continuous model considered in \cite{sarkar}. Here points of a Poisson point process are connected to their $k\ge1$ nearest neighbors. To study this percolation process, the authors develop a rigorous finite-size criterion, which is then verified ``with high probability'' through a simulation. 

We next aim to transfer this finite-size approach to our situation, similar adaptation were considered in \cite{stepDisk,stacey}.
The key idea is to identify edges $e$, where it is already possible to say after a finite number of rounds whether $e \in \mc E_\ff$. More precisely, an edge is after $n \ge 1$ rounds
\been
\im  \emph{certainly vacant}  if it is reinforced at most $n - 1$ times;
\im \emph{potentially occupied} if it is reinforced precisely $n$ times;
\im \emph{certainly occupied}  if it is reinforced at least $n + 1$ times. Moreover, if a node is incident to two certainly vacant edges, then the remaining two edges are also \emph{$n$-certainly occupied}.
\enen
The reason for this terminology is that as shown in the proof of Lemma \ref{lem:stdom}, if an edge is certainly occupied, then it is eventually reinforced  in every round. Similarly, if an edge is certainly vacant, it is only reinforced a finite number of times.

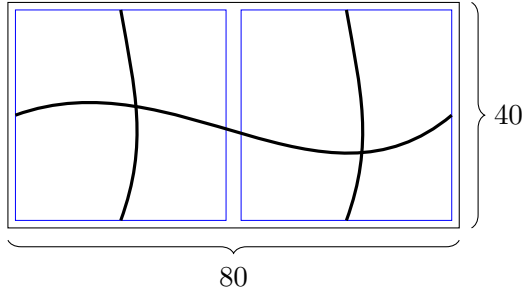
\begin{figure}[t]
\begin{tikzpicture}
    % Rectangle
    \draw (0, 0) rectangle (6, 3);

    % two smaller rectangles in blue
    \draw[blue, very thin] (0.1, 0.1) rectangle (2.9,2.9);
    \draw[blue, very thin] (3.1, 0.1) rectangle (5.9,2.9);
    
    % Horizontal crossing (curved)
    \draw[very thick] (0.1, 1.5) to[out=20, in=-140] (5.9, 1.5);

    % Vertical crossings (curved)
    \draw[very thick] (1.5, 0.1) to[out=70, in=280] (1.5, 2.9);
    \draw[very thick] (4.5, 0.1) to[out=70, in=280] (4.5, 2.9);
    \draw [decorate,decoration={brace,amplitude=5pt,mirror,raise=4ex}]
  (0,0.5) -- (6,0.5) node[midway,yshift=-3em]{80};
      \draw [decorate,decoration={brace,amplitude=5pt,mirror,raise=4ex}]
  (5.5,0) -- (5.5,3) node[midway,xshift=3em]{40};
\end{tikzpicture}
	\caption{Illustration of the crossings in an $n$-open rectangle}
	\label{fig:cross}
\end{figure}

Now,  we say that the rectangle $80 \times 40$ is \emph{$n$-open} if there exist three types of crossings of the central $(80 - 2n) \ti (40 -2n)$ rectangle with edges that are certainly occupied after $n$ rounds:
\been
\im a horizontal crossing of the $((80 - 2n) \ti (40 -2n))$-rectangle;
\im vertical crossings of the left and right $((40-2n) \ti (40 - 2n))$-squares in the central rectangle.
\enen
Figure \ref{fig:cross} shows the crossings of the central rectangle. There is no particular conceptual reason for the values 80 and 40. However, we found them to be convenient for the simulations in Section \ref{sec:sim}. Henceforth, we often think of the $(80 \times 40)$-rectangle as being a horizontal edge in a coarse-grained lattice.

Having introduced the notion of $n$-openness, we define $\eff(n)$ as the event that the coarse-grained edges $(0,0) \to (1, 0)$ is $n$-open. Then, we show that if $\P(\eff(n)) > 0.8457$, then $\P(\EE_{\ff} \text{ percolates}) = 1$. The deeper reason behind the value 0.8457 is the recent paper \cite[Theorem 1]{balister}, which shows that any 1-dependent 2D bond percolation model percolates if the marginal probability exceeds this value.

%
%THM A=FF
%
\bet[Finite-size criterion]
\label{thm:reinf}
Let $d= k = 2$ and  $n\ge 1$ be such that $\P(\eff(n)) > 0.8457$. Then,
$$\P(\EE_{\ff} \text{ percolates}) = 1.$$ 
\ent

Using Monte Carlo simulation,  we show that with a certainty exceeding $1 - 10^{-300}$, the finite-size crossing probability in Theorem \ref{thm:reinf} indeed exceeds the threshold 0.8457 of 1-dependent percolation.

%
%FIN ALPHA
%

\subsection{The case $\a < \ff$}
\label{ssec:fin}

 Next, we deal with the case $\alpha < \ff$. We go through the algorithm for the weight evolution similarly as in the case where $\a = \ff$ and again initialize the weights with $W_e(0) = 1$ for every $e \in E$. Then, for each $t\ge 1$, we carry out the following updates
\been
\im At every site $z \in \Z^d$, we select one tuple $\s \in E_z^{(k)}$ whose edge weights are increased by 1.
\im These edges are selected with probability proportional to $\prod_{e \in \s} W_{e}(t - 1)^\a$. 
\enen 
In the limiting case $\a = \ff$, we recover the simplified reinforcement mechanism discussed above. Therefore, it is plausible to extend Conjecture \ref{conj:reinf} from $\a = \ff$ to the case of large but finite $\a$. 

%
%CONJ REINFA
%
\begin{conjecture}[Percolation of reinforced 2-out percolation on $\Z^2$; finite $\alpha$]
\label{conj:reinfa}
Assume that $d= k = 2$.  Then, there exists $\a_0 > 0$ such that  $\P(\EE_{\ff} \text{ percolates}) = 1$ whenever $\a > \a_0$.
\end{conjecture}

As before, we add independent sprinkling on the edges with a sprinkling probability $p_{\ms{spr}} > 0$. As in the case of $\a = \ff$, we consider then the union $\EE_{\ff, \e} := \EE_\ff \cup \EE_{\e, \ms{spr}}$ 
of edges that are reinforced infinitely often together with the sprinkled edges.

%
%CORR
%
Ideally, as a first step, we would like to present a finite-size criterion as in Theorem \ref{thm:reinf}.  The key difficulty is that now, it is more complicated to identify edges where we know that they will be reinforced either finitely or infinitely often. For instance, if an edge has weight 1 after 10 rounds, it could still happen (although extremely unlikely) that it will still be reinforced in each of the following 100 rounds. 

Hence, we need to introduce a process of corrupted vertices to take into account such effects. More precisely, let $p_* \in[0, 1]$ be a fixed (small) corruption probability. Then, each vertex incident to at least one edge of weight at most $n$ after $n$ rounds is declared \emph{$n$-corrupted} independently with probability $p_*$.
Having introduced the corrupted vertices, we can now extend the concept of certainly occupied edges. More precisely, after $n \ge 1$ rounds, we say that an edge that is not incident to a corrupted vertex is
\beit
\im \emph{certainly vacant} if it is reinforced at most $n - 1$ times;
\im \emph{certainly occupied} if it is not certainly vacant and is incident to a node which itself is incident to two certainly vacant edges. 
\enit

To introduce the finite-size criterion, we now proceed in the same way as in $\a = \ff$. Again, to avoid boundary effects, we consider a central rectangle $(80 - 2n) \ti (40 -2n)$ rectangle with edges that are certainly occupied after $n$ rounds. Then, we say again that it is \emph{$n$-open} if there exist the three types of crossings considered in the case $\a = \ff$. Now, we say that a horizontal edge between $(x, y)$ and $(x + 1, y)$ is \emph{$n$-open} if the rectangle $[40x-40, 40x + 40) \ti [40y-20, 40y+20)$ has this property in the above sense. Similarly, we define the openness of vertical edges.  Now, we let $\efa = \efa(n)$ denote the event that the coarse-grained edges $(0,0) \to (1, 0)$ is $n$-open with corruption probability $p_*$ set as
\begin{align}
	\label{eq:psa}
p_*(\alpha, n) := 1\wedge\Big( (n - 1)^\a\sum_{j \ge n }{j^{-\a}}\Big).
\end{align}
We note that the infinite sum cannot be evaluated in closed forms. However, it would be possible to present closed-form  bounds through suitable integral bounds. {For fixed $n$, we can (and will) choose $\a$ so that $p_*$ is sufficiently small.}

%
%THM A<FF
%
\bet[Percolation of reinforced 2-out percolation on $\Z^2; \a \gg 0$]
\label{thm:reinf1}
Let $d= k = 2$.  Furthermore, assume that $\alpha > 0$ and $n$ are such that  $\P(\efa(n)) > 0.8457$. Then, $\P(|\mc E_\a| = \ff ) = 1$.

\ent

We note that the bound on $p_*$ is the result of a comparison with a strongly-reinforced P\'olya urn process, which we elaborate on in Lemma \ref{lem:cor} below. Again, in Section \ref{sec:sim}, we provide overwhelming evidence from Monte Carlo simulations that $\P(\efap(n)) > 0.8457$.

Finally, we stress again that the difficulty of analyzing the reinforced model is that it exhibits long-range dependencies and does not satisfy the FKG inequality. This means that the vast majority of the standard percolation arguments break down for this model. Nevertheless, in the case $\a = \ff$, the model enjoys the following intriguing domination property. 

%
%THM WRINF
%
\bepr[Domination of $\mc E_\ff$ when $\a = \ff$]
\label{pr:dom}
If $\a=\ff$, then the random set $\mc E_\ff$ stochastically dominates $E \setminus \mc E_\ff$. 
\enpr

It is unclear if the stochastic-domination property also holds for $\a < \ff$. However, we believe that an analog of Proposition \ref{pr:dom}  could be possible for a clever choice of finite reinforcement function.

%% file: ind.tex
\section{Proof of Proposition \ref{pr:dom}}
\label{sec:dom}

We will describe a new colored process, equivalent to our original reinforcement process for $\alpha=\infty$, from which the derivation of Proposition \ref{pr:dom} will be straightforward.
In this colored process, each edge $e\in E(\Z^d)$ has a blue counter $b_t(e)$ and a red counter $r_t(e)$ that increase with time. At the beginning we set $r_0(e)=b_0(e)=1$ for all edges. 

At each step, every vertex chooses $d$ of its adjacent edges, which increase their blue counter by $1$, and the other two increase their red counter by $1$. The choice is made according to the same reinforcement mechanism as in our regular process --- we choose the edges incident to $v$ with the highest $b(e)$ value to increase theis blue counter while all the other $d$ incident edges increase their red counter, breaking ties by choosing uniformly at random. Set 
$$B_n:=\{e: b_n(e)> n+1\},\; R_n(e):=\{e: r_n(e)> n+1\},\;\text{ and }U_n:=\{e: b_n(e)= n+1\}.$$ 
Finally set 
$$B_\infty :=\bigcup_n B_n,\; R_\infty :=\bigcup_n R_n\; \text{and } U_\infty=\bigcap U_n.$$
We make the following simple observations:
\bel
\label{lem:stdom}

\begin{enumerate}
    \item For every edge $e$ and every $n\geq 0$,  $b_n(e)+r_n(e)=2n+2$.
    \item  For every $n \ge 0$, it holds that $B_n\sqcup R_n \sqcup U_n =\Z^d$, where $\sqcup$ means that the union is disjoint.
    \item $B_n,R_n$ are increasing sequences, $U_n$ is decreasing. 
    \item $ B_\infty \cup R_\infty\cup U_\infty =\Z^d $.
    \item $B_n$ and $R_n$ have the same distribution, $B_\infty$ and $R_\infty$ have the same distribution.
    \item $B_\infty \cup U_\infty = \{e: b_n(e)\to \infty\}$, $R_\infty \cup U_\infty = \{e: r_n(e)\to \infty\}$.
\end{enumerate}
\enl

Before proving Lemma \ref{lem:stdom}, we explain how to deduce Proposition \ref{pr:dom}, where the main step will be to identify the models in such a way that $\EE_\ff = B_\ff \cup U_\ff$. This identification also gives that $R_1 \cap \EE_\ff = \es$, which implies that our model is stochastically dominated by the independent $k$-out model with $k=d$.

\begin{proof}[Proof of Proposition \ref{pr:dom}]
The connection to our reinforcement model comes by looking only at the blue counters. Under this identification we get that $B_n$ is the set of all certainly occupied edges and $R_n$ is the set of certainly vacant edges at time $n$. Also $\EE_{\ff} = B_\infty \cup U_\infty$.
The claim now follows directly from clause $(4)$ of Lemma \ref{lem:stdom}.
\end{proof}

Now, we give the proof of Lemma \ref{lem:stdom}.
\begin{proof}[Proof of Lemma \ref{lem:stdom}]
\begin{enumerate}
    \item The assertion follows from the fact that each edge gets either a red counter or a blue counter from each of its endpoints.
    \item Clause $(1)$ implies that $R_n,B_n,U_n$ are a disjoint partition of $\E(\Z^d)$.
    \item We will show this for $B_n$, and $R_n$ follows by symmetry. Take some edge $e=(v_+,v_-)\in B_n$ and consider the first time $m< n$ for which $e\in B_{m+1}$. Then at time $m$, the edge $e$ got a blue counter from both its endpoints. Therefore, $b_m(e)$ was among the $d$ largest blue counters of edges incident to $v_+$ (and similarly for $v_-$).  
    Moreover, the blue counter of any edge can increase by at most $2$ at every step.
    Therefore, $b_{m+1}(e)$ is strictly larger than the $d$ lowest blue counters of edges incident to $v_+$. As a consequence,  the edge $e$ will keep being chosen by both  $v_+$ and $v_-$ for all rounds $k \ge m + 1$. This implies that $e\in B_k$ for all $k\ge  m+ 1$. Finally, the fact that $U_n$ is decreasing now follows from clause $(2)$.
     \item  Clause $(3)$ implies that the unions defining $B_\infty$ and $R_\infty$ are increasing, and that the intersection defining $U_\infty$ is decreasing. Therefore, the statement follows from clause $(2)$.
    \item By clause $(1)$, choosing the edges with the highest blue counter is the same as not choosing those with the highest red counter. This makes the definition of the process symmetric w.r.t. flipping the colors. 
    \item This follows from a similar argument to the monotonicity clause $(2)$. 
\end{enumerate}
\end{proof}

%% file: osss.tex
%
%SEC OSSS
%
\section{Proof of Theorem \ref{thm:osss}}
\label{sec:osss}

In this section, we prove Theorem \ref{thm:osss}. That is, we show that percolation occurs in the dual or original lattice after adding arbitrarily sparse independent sprinkling. In \eqref{eq:eeff}, {we defined the set $\EE_\ff$ of edges that are eventually reinforced in every round}. One challenge of dealing with this set is that it has an unbounded range of dependencies. To overcome this difficulty, we now introduce an approximated version $\EE_N$ of the set $\EE_\ff$ after a finite number of $N \ge 1$ rounds. 
 We say that an edge $e$ is \emph{$N$-potentially occupied} (and write $e\in\EE_N$) given the configuration after $N$ rounds if it is still possible with positive probability that $e \in \EE_\ff$. 
Further, we call an edge $e$ \emph{strictly $N$-potentially occupied} if  $ e\in \EE_N \sm \EE_\ff$.  The main idea is to show percolation for each of the approximated sets $\EE_N$ and then to use a stochastic domination result to control the difference $ \EE_N \sm \EE_\ff$.

%
%DEMBIN
%
This strategy is similar to that used in \cite{dembin}, where the authors prove a similar result for Voronoi percolation. The key similarity is that both Voronoi percolation and our model are spatial percolation models with exponential decay of correlations.   

To summarize, the proof of Theorem \ref{thm:osss} relies on two central results, namely Propositions \ref{pr:ss} and \ref{pr:comp} below.

%
%PRF SS
%
\bepr[Sharp threshold for approximation]
\label{pr:ss}
Fix $N \ge 1$ and let $\ps( N)$ denote the critical sprinkling probability in the $N$-approximated model. Then, for $p < \ps( N)$, the diameter of the clusters has exponentially decaying tails. 
\enpr

Let $S_N$ denote the collection of sites that are incident to at least one edge that is strictly $N$-potentially occupied.
%
%PRF COMP
%
\bepr[Stochastic domination]
\label{pr:comp}
Let $N \ge 2$. Then, $S_N$ is dominated by a Bernoulli site percolation process with marginal probability $\de_N = 6^{-N/2}$. 
\enpr

We first discuss how to prove Theorem \ref{thm:osss} subject to the validity of Propositions \ref{pr:ss} and \ref{pr:comp}. Then, we prove these propositions in the subsequent subsections.

%
%PRF OSSS	
%
\bep[Proof of Theorem \ref{thm:osss}]
{Proposition \ref{pr:comp}} gives a stochastic domination of $S_N$ by a Bernoulli site percolation process with marginal probability $\de_N$. We first claim that from this result, we can conclude that the edges incident to sites in $S_N$ are stochastically dominated by a Bernoulli bond percolation process with probability $\sqrt{\de_N}$. Indeed, consider an arbitrary finite edge set $F \su E$. Let $V_F$ denote the set of left or lower endpoints of the edges in $F$. Then, by Proposition \ref{pr:comp},
$$\P(F \su S_N) \le \de_N^{|V_F|} \le \de_N^{|F|/2},$$
thereby proving the claimed stochastic domination. 
A consequence is that $\EE_{N(\e)} \su \EE_{\ff} \cup \EE_{\e, \ms{spr}}$ provided that $N(\e)$ is so large that $\de_N^{1/2} \le \e$. {We also note that the independent superposition of two Bernoulli bond percolation processes with parameter $\e >0$ is again a Bernoulli bond percolation process, whose parameter is then $1 - (1 - \e)(1 - \e)  = 2\e -\e^2$.}

{
Now, to derive a contradiction, assume that for some $\e > 0$
$$\P(\EE_\ff \cup \EE_{2\e - \e^2, \ms{spr}}\text{ percolates}) = \P(\EE_\ff \cup \EE_{2\e^2- \e, \ms{spr}}\text{ percolates dually}) = 0.$$ 
Then, Proposition \ref{pr:comp} implies that,
$$
\P(\EE_{N(\e)} \cup \EE_{\e, \ms{spr}} \text{ percolates})  = \P(\EE_{N(\e)} \cup \EE_{\e, \ms{spr}}\text{ percolates dually}) = 0. 
$$
In other words, $\e$ is a subcritical sprinkling probability both with respect to original and dual connectivities of the $N(\e)$-approximated model. In particular,  $\EE_{N(\e)} = \EE_{N(\e)} \cup \EE_{0, \ms{spr}}$ is in the strictly subcritical regime, where Proposition \ref{pr:ss} guarantees the exponential tail decay of the cluster sizes.
}
%$\P(\EE_{N(\e)} \text{ percolates weakly}) = 0$. 
In particular, by Proposition \ref{pr:ss}, both the sizes of the clusters with respect to the original connectivities and with respect to the dual connectivities  have exponentially decaying tails. 
Consequently, the probability to have left-right vacant crossings of large squares tends to 1 as the side length of the square tends to infinity. 

Now, by Proposition \ref{pr:dom}, the occupied edges dominate the vacant ones in $\EE$, thereby contradicting the exponential decay of the cluster sizes. 
\enp

%
%FACTOR ENCODING
%
We will prove Propositions \ref{pr:ss} and \ref{pr:comp} in Sections \ref{ss:ss} and \ref{ss:comp} below. For both results, it is convenient to describe more precisely how our model can be constructed as a factor of a model with independent inputs. 
First, we introduce an iid sequence $\{U_e^{\ps}\}_{e \in E}$ that encodes the sprinkling of edges. More precisely, each $U_e^{\ps}$ is uniformly distributed in $[0,1]$ and the set of sprinkled edges is then given by 

% there is additional randomness used for sprinkling the edges. To incorporate it into the model formally, we introduce an iid sequence $\{U_e^{\ps}\}_{e \in E}$ of uniformly distributed random variables on $[0,1]$ such that the set of sprinkled edges is given by
$$X := \{e \in E\co U_e^{\sp} \le \ps\},$$
   Second, we attach to each site $z \in \Z^2$ a sequence of iid decision variables $U_z^{(1)}, U_z^{(2)}, \dots$. All variables are independent for different values of $z$ and also independent of the sprinkling variables $X$.
   
   Given the configuration of edge reinforcements up to iteration $n$, the variable $U_z^{(n)}$ is used to decide which of the edges are reinforced by $z$ in iteration $n$. We write 
    $$Y_z := \{U_z^{(k)}\}_{k\ge 1}$$ 
    for the sequence of all decision variables at $z$. Then, we write $Y := \{Y_z\}_{z \in \Z^2}$ for the collection of all $Y_z$.   Hence, the entire randomness of the model can be encoded in the pair 
    $$Z:=(X, Y).$$

    %
    %SS COMP
    %
    \subsection{Proof of sharp thresholds for approximation -- Proposition \ref{pr:ss}}
    \label{ss:ss}
    We will deduce Proposition \ref{pr:ss} from a differential inequality in the spirit of \cite[Lemma 1.7]{voronoi}. We also write $\th_n(\ps(N))$ for the probability that in the sprinkled model the connected percolation component of the origin has $\ell_\ff$-diameter exceeding $n$, {which we denote as $o\lrsa \pa B_n $}.
    
    %
    %DIF LEM
    %
    \bepr[Differential inequality]
    \label{sup_crit_lem}
    There exists $c_{\ms{Diff}} > 0$ such that for every $n \ge 1$  and $\ps > 0$ we have 
    $$\f{\d}{\d \ps}\th_n(\ps) \ge c_{\ms{Diff}} \f n{\sum_{s \le n}\th_s(\ps)}\th_n(\ps) (1 - \th_n(\ps)).$$
    \enpr
    Applying the OSSS  inequality (after O'Donnell, Saks, Schramm and Servedio, \cite{osss}) to an algorithm $T$ determining $f_n :=\one\{o\lrsa \pa B_n\} $ gives that
    \begin{equation}
            \label{osss}
    	\th_n(\ps) (1 - \th_n(\ps)) = \Var(f_n(Z)) \le \sum_{z \in \Z^d} \de_z(T) \Inf_z^Y + \sum_{e \in E(\Z^d)}  \de_e(T) \Inf_e^X ,
    \end{equation}
    where
    \been
    \im[(i)] $\de_z(T) := \P(T\text{ reveals }Z_z)$ is the probability that the algorithm $T$ reveals the value of $X_z$, and
    \im[(ii)] $\Inf_z^Y := \P\big(o\lrsa \pa B_n \ne o\lrsa^z \pa B_n\big)$ denotes the \emph{(resampling) influence}, where $\lrsa^z$ is the percolation event after replacing the randomness $Y_z$ at $z$ by an independent copy, while leaving $X_z$ fixed.
    \enen
    The quantities $\de_e(T)$ and $\Inf_e^X$ are defined similarly. Henceforth, we rely on a specific algorithm $T$ from \cite{voronoi}, which involves an additional randomization through a uniform integer $U \in \{1, \dots, n\}$.  For the convenience of the reader, and to make the manuscript self-contained, we briefly recall the idea behind this algorithm. The idea is to explore the clusters by starting from $\pa B_U$. More precisely, we proceed as follows:
    \been
   \im Reveal the value of $U_e^{\ps}$ for  all edges incident to $\pa B_U$. Also reveal $Y_z$ for all $z \in \Z^2$ at $\ell_\ff$-distance at most $N$ from a point in $\pa B_U$.
    \im Suppose that the values of $X_{e_1}, \dots, X_{e_k}$ and of $Y_{z_1}, \dots, Y_{z_\ell}$ have already been revealed at the start of iteration $t$ of the algorithm. Let $\mc C_t$ be the union of all connected components that are revealed by this iteration. Then, according to some arbitrary  rule, we pick out some unrevealed site $z_0$ that is incident to an edge in $\mc C_t$. We then reveal the values of $X_e$ for any edge adjacent to $z_0$. We also reveal the values of $Y_z$ for all $z \in \Z^2$ at $\ell_\ff$-distance at most $N$ from $z_0$.
    \enen

    The proof of Proposition \ref{sup_crit_lem} relies on the following two central auxiliary results (Lemmas \ref{lem:inf} and \ref{lem:rev} below) concerning a comparison of the two influences, and a bound on the revealment probabilities, respectively. 
    
    %
    %LEM INF
    %
    \bel[Comparison of influences]
    \label{lem:inf} Let $N \ge 1$. Then, there is $c_1 = c_1(d, N) > 0$ such that 
    for every  $z \in \Z^2$,
    $\Inf_z^Y \le c_1 \sum_{e \su B_N(z)}\Inf_e^X.$
    \enl

    %
    %LEM REV
    %
    \bel[Bound on revealment probabilities]
    \label{lem:rev}
    Let $N \ge 1$. Then, there is $c_2 = c_2(d, N) > 0$ with the following property.
    Let $T$ the randomized algorithm starting the exploration from $\pa B_U$, with $U$ uniform in $\{1, \dots, n\}$. Then,
    $$\big(\sup_{e \in E}\de_e(T)\big) \vee \big(\sup_{z \in \Z^2}\de_z(T)\big) \le \f{c_2\sum_{s \le n}\th_s(\ps)}n$$
    \enl
    
    %
    %PRF SUP CRIT
    %
    Before proving Lemmas \ref{lem:inf} and \ref{lem:rev}, we explain how to conclude the proof of Proposition \ref{sup_crit_lem}.
    \bep[Proof of Proposition \ref{sup_crit_lem}]
    First, by the OSSS inequality \eqref{osss} and Lemmas \ref{lem:inf} and \ref{lem:rev}, we obtain that
    $$	\th_n(\ps) (1 - \th_n(\ps)) \le c_1c_2|B_N(o)|n^{-1}{\sum_{s \le n}\th_s(\ps)}\sum_{e \in E(\Z^d)} \Inf_e^X.$$
    Now, Russo's formula gives that 
    $$\sum_{e \in E(\Z^d)} \Inf_e^X = \f{\d}{\d \ps}\th_n(\ps),$$
    thereby concluding the proof.
    \enp
    
    %
    %PRF INF
    %
    Hence, we now prove the auxiliary results Lemmas \ref{lem:inf} and \ref{lem:rev}, starting with Lemma \ref{lem:inf}. The idea is similar to \cite{dembin}, where a comparison between different forms of pivotality also plays a crucial role. However, our situation is a bit simpler in the sense that $N$-approximated model has the bounded dependence range $N$.
    \bep[Proof of Lemma \ref{lem:inf}]
    For $\s \in \{0, 1\}$, we introduce the event 
    $$E_{\ms{coarse}, \s} := \{f_n(Z_\s)  = \s\},$$
    where $Z_\s$ denotes the configuration obtained from $Z$ by setting $1 - U_e^{\sp} := \s$ for every $e \in E \cap B_N(z)$. 

    Now, using that the $N$-approximated model has dependence range $N$, we see that if the resampling of $Y$ at $z$ changes the value of $f_n$, then $E_{\ms{coarse}} := E_{\ms{coarse}, 0}\cap E_{\ms{coarse}, 1}$ must occur. In other words, $\Inf_z^Y \le \P(E_{\ms{coarse}})$. Hence, it suffices to relate the probability of $E_{\ms{coarse}}$ with the sum of the influences $\Inf_e^X$ for $e \in B_N(z)$.
    
    To achieve this goal, we write $X^*$ for the configuration of the model with the sprinkled edges $X_e$ replaced by independent copies for every $e \in B_N(z)$. Then, we let $E_{\ms{fine}, -}$ denote the event that in $X$ none of the edges in $B_N(z)$ is sprinkled. We also let $E_{\ms{fine}, +}$ denote the event that in the resampled configuration all edges in $B_N(z)$ are sprinkled. Then, by independence and the definition of $E_{\ms{coarse}}$, we have
    $$c\P(E_{\ms{coarse}}) \le \P\big(E_{\ms{coarse}} \cap E_{\ms{fine}, -} \cap E_{\ms{fine}, +}\big) \le \P\big(f_n(X, Y) \ne f_n(X^*, Y)\big)$$
    for some $c = c(N)>0$.     Noting that the right-hand side is bounded above by $\sum_{e \su B_N(z)}\Inf_e^X$ concludes the proof.
    \enp

    %
%PRF REV
%
Hence, it remains to establish the revealment bounds asserted in Lemma \ref{lem:rev}. Here, we proceed similarly as in \cite{voronoi}, noting that again the finite-range property of our model simplifies the argument. 
\bep[Proof of Lemma \ref{lem:rev}]
We only bound the revealment probability $\de_z(T)$, as the bound on $\de_e(T)$ is easier. The key observation is that if the randomized algorithm starts from exploring $\pa B_U$ with $U= m$ for some $m \le n$, then the following holds. If $T$ reveals the state of $Y_z$ for some site $z\in \Z^2$, then $\pa B_m \lrsa B_N(z)$. Therefore,
$$\P(T \text{ reveals }z) \le \P\big(\pa B_m \lrsa B_N(z)\big) \le \th_{|m - |z|_\ff| - N},$$
where we use the convention $\th_k = 1$ for $k \le 0$.
Thus, picking $m \in \{1, \dots, n\}$ uniformly at random, we obtain that
$$\de_z(T) \le \f 1n\Big(2N + \sum_{m \le n} \th_{m}\Big) \le \f cn\sum_{m \le n} \th_{m},$$
for a suitable $c = c(N) > 0$. This concludes the proof.
\enp

    %
    %SS COMP
    %
    \subsection{Proof of stochastic  domination of approximation -- Proposition \ref{pr:comp}}
\label{ss:comp}
To prove Proposition \ref{pr:comp}, we need to describe more precisely how the decision variables $\{U_z^{(k)}\}_{k \ge 1}$ determine the state of the edges.

%
%PRF COMP
%
\bep[Proof of Proposition \ref{pr:comp}]
 Our construction depends on the parity of $k$, and we first discuss the case of odd $k$. Henceforth, we implicitly consider $\Z^2$ to be endowed with a checkerboard pattern so that we can speak of black and white sites. 

For a black site $z \in \Z^2$, we let $U_z^{(k)}$ encode in some arbitrary way a random selection of two of the highest-weight edges incident to $z$. Next, to determine the reinforcements at a white site $z' \in \Z^2$, we work conditioned on the edge weights up to iteration $k$ and also on the decision variables $U_{z}^{(k)}$ for all black sites $z \in \Z^2$.  We note that the probability space $\Omega$ of selecting two highest-weight edges incident to $z'$ consists of $\binom42 = 6$ elements. 

{First, consider the case where, in step $k$, at least 2 of the black neighbors of $z'$ reinforced  their corresponding edge to the white vertex $z'$. Then, there is a probability of at least $\de = 1/6$ that, in step $k$, the vertex $z'$ selects the same two edges, which are therefore reinforced infinitely often, see clause (3) in Lemma \ref{lem:stdom}. 

Second, consider the case when there are at least 2 neighboring black vertices of $z'$ that did not reinforce their corresponding edge to the white vertex $z'$. Then, again with probability of at least $\de = 1/6$, in step $k$, the vertex $z'$ selects the other two edges. Then, again by Lemma~\ref{lem:stdom}, these two edges are never reinforced again. Therefore, the  vertex $z'$ is no longer incident to any strictly $N$-potentially occupied edges.

 Writing $Y_-^{(k)}$ for the collection of decision variables of all sites up to step $k - 1$ together with the decision variables at the black sites at step $k$, we conclude that there exists a state $\om_*(Y_-^{(k)}) \in \Omega$ such that almost surely, i) $\P(U_{z'}^{(k)} = \om_*(Y_-^{(k)})\ba Y_-^{(k)}) \ge \de$ and ii) if $U_{z'}^{(k)} = \om_*(Y_-^{(k)})$, then after step $k$, all potentially reinforced edges incident to $z'$ are in fact reinforced infinitely often.
}

Taking into account these observations, we now present a two-step construction of the variables $U_{z'}^{(k)}$. For this, let $\{W_z^{(k)}\}_{k \ge 1, z \in \Z^2}$ be an iid sequence of Bernoulli random variables with parameter $\de$. If $W_z^{(k)} = 1$, then we let $U_z^{(k)}$ be the state $\om_*(Y_-^{(k)})$. Otherwise, if $W_z^{(k)} = 0$, then we let $U_z^{(k)}$ be the state $\om_*(Y_-^{(k)})$ with probability $q_0(Y_-^{(k)}) := \big(\P(U_z^{(k)} = \om_*(Y_-^{(k)})\ba Y_-^{(k)}) - \de\big)/(1- \de)$, and with probability $1 - q_0(Y_-^{(k)})$ the random variable $U_z^{(k)}$ is sampled according to the conditional distribution $\mc L\big(U_z^{(k)} | U_z^{(k)} \ne \om_*(Y_-^{(k)}), Y_-^{(k)}\big)$.

This describes the construction for odd steps $k$. For even $k$, we proceed in precisely the same manner except that the roles of black and white sites are interchanged. In particular, if an edge $e$ incident to a black site $z$ is contained in $\EE\sm\EE_N$, then $W_z^{(1)} = W_z^{(3)} = \cdots = W_z^{(N)} = 0$. This proves the asserted domination.
\enp

%% file: reinf.tex
%
%SEC REINF
%
\section{Proofs of Theorems \ref{thm:reinf} and \ref{thm:reinf1}}
\label{sec:reinf}

%
%A = FF
%
We  deal separately with the cases $\a = \ff$ and $\a < \ff$. We start with $\alpha=\ff$, where we use a result for 1-dependent percolation with sufficiently high marginal probabilities from \cite{balister}. We apply this result to the coarse-grained model. The key observation is that if the coarse-grained edge $(0, 0) \to (1, 0)$ is $n$-open, then we are guaranteed three crossings of infinitely-reinforced edges in the rectangle $[0, 80) \ti [0, 40)$. Namely, 
\been
\im a horizontal crossing of the $((80 - 2n) \ti (40 -2n))$-rectangle;
\im vertical crossings of the left and of the right $((40-2n) \ti (40 - 2n))$-squares inside the central rectangle.
\enen
Moreover, the existence of such crossings depends only on coarse-grained edges sharing at least one of the end points. Hence, we conclude from planarity that any path of  $n$-open coarse grained edges gives rise to a path of infinitely-reinforced edges. Hence, it suffices to establish the percolation of the coarse-grained model, which we do now.

%
%PRF A= FF
%
\bep[Proof of Theorem \ref{thm:reinf}; $\a = \ff$]
The collection of $n$-open edges defines a 1-dependent family. Now, by our assumption on the finite-size criterion, a coarse-grained edge is open with marginal probability exceeding 0.8457. Now, we can conclude the proof by invoking \cite[Theorem 1]{balister}, which states that any 1-dependent site percolation model with marginal probability exceeding 0.8457 percolates.
\enp

%
%PRF A < FF
%
We now turn to the case $\a < \ff$. We first argue that percolation of the $n$-open edges in the coarse-grained model implies percolation of the infinitely-reinforced model in the original model.

%
%LEM COR
%
\bel[Stochastic domination property of certainly occupied edges]
\label{lem:cor}
Let $n, \a > 1$ and  assume 
$$p_* > (n - 1)^\a\sum_{j \ge n }{j^{-\a}}.$$
Then, the process of edges that are reinforced only finitely often in the original model is stochastically dominated by the process of edges that are not certainly occupied.
\enl

%
%PRF COR
%
\bep
The idea of the proof is as follows. If after $n$ rounds an edge $e$ has weight at most $n$, then for large $\a$, it is highly likely never to be chosen again. Then, loosely speaking,  the vertex corruption is used to account for the possibility of this exceptional event.

To make this precise, we let $E_e(n)$ denote the event that the edge $e$ with weight at most $n - 1$ in round $n$ is reinforced by at least one of its incident vertices in some round $j \ge n$. Then, we claim that $(\one\{E_e(n)\})_e$ is stochastically dominated by a directed Bernoulli bond percolation process with probability $p_*$. To achieve this goal, we note that at the beginning of  round $j \ge n$ there are at least 2 high-weight edges, i.e., edges of weight at least $j$. Hence, the probability that one of the low-weight edges $e$ is reinforced from one of its incident vertices is at most 
$$2\f{(n - 1)^\a}{2j^\a} = \Big(\f{n - 1}j\Big)^\a.$$
We stress that this upper bound holds irrespective of the weight evolution at any of the other directed edges. Hence, the union bound shows that, as asserted, the edge process $(\one\{E_e(n)\})_{e \in E}$ is dominated by a Bernoulli bond process with probability $p_*$.
\enp

%
%PRF A < FF
%
\bep[Proof of Theorem \ref{thm:reinf}; $\alpha <\ff$]
After the reduction step from Lemma \ref{lem:cor}, the proof is very similar to that in the case $\a = \ff$. Arguing as in the case $\a = \ff$, Lemma \ref{lem:cor} guarantees that  any path of  $n$-open coarse grained edges gives rise to a path of infinitely-reinforced edges. Moreover, we conclude from the condition \eqref{eq:psa} that the probability for a coarse-grained edge to be $n$-open exceeds 0.8457. The theorem then follows, as for $\a=\ff$, by invoking \cite[Theorem 1]{balister}.
\enp

%% file: sim.tex
%
%SEC REINF
%
\section{Numerical evidence for the finite-size criteria}
\label{sec:sim}
Now, we give numerical evidence that the finite-size criteria from Theorems \ref{thm:reinf} and \ref{thm:reinf1} are satisfied.

\begin{figure}[!htpb]
        \centering                                       
\includegraphics[width=0.48\textwidth]{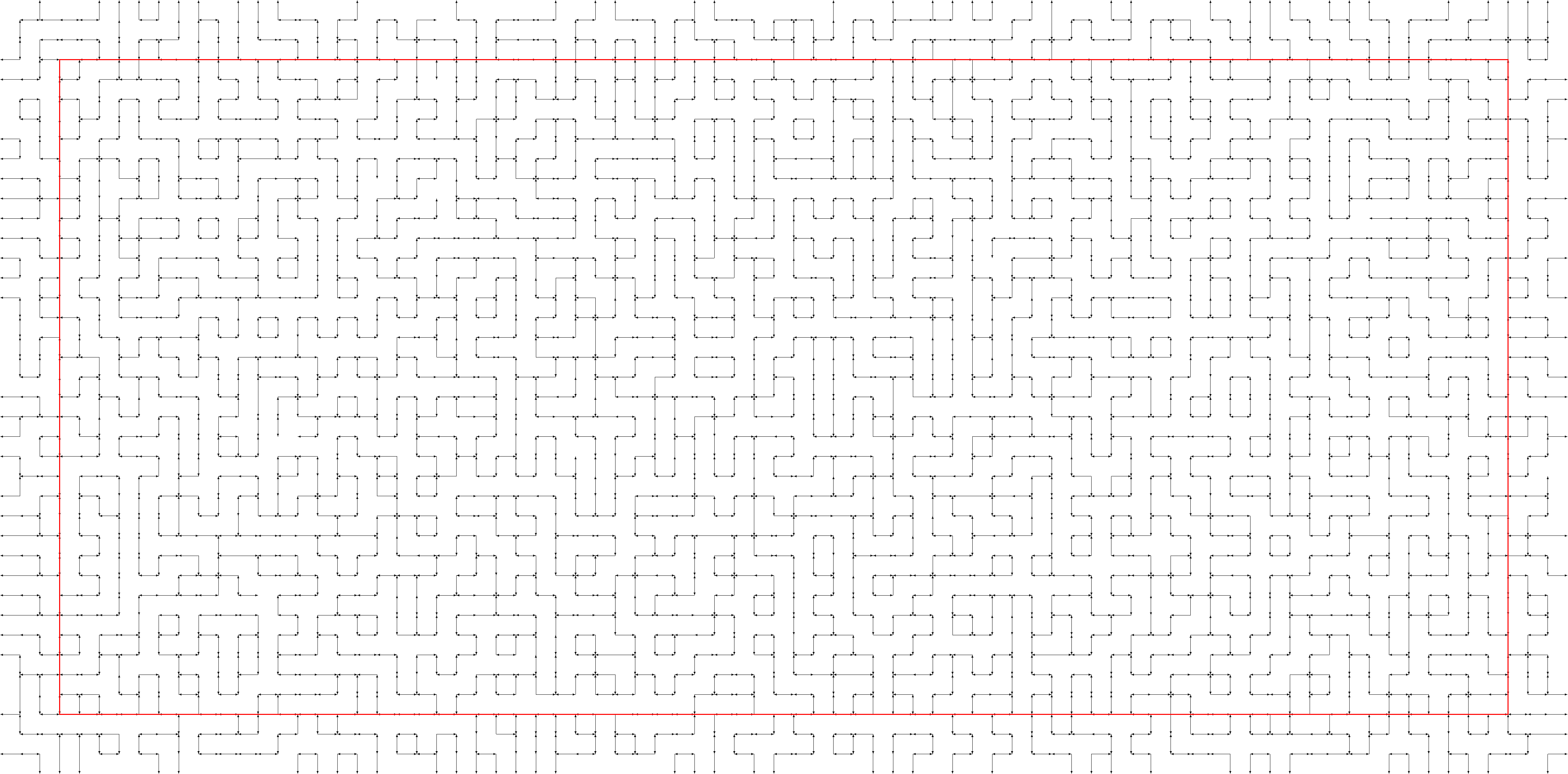}
\includegraphics[width=0.48\textwidth]{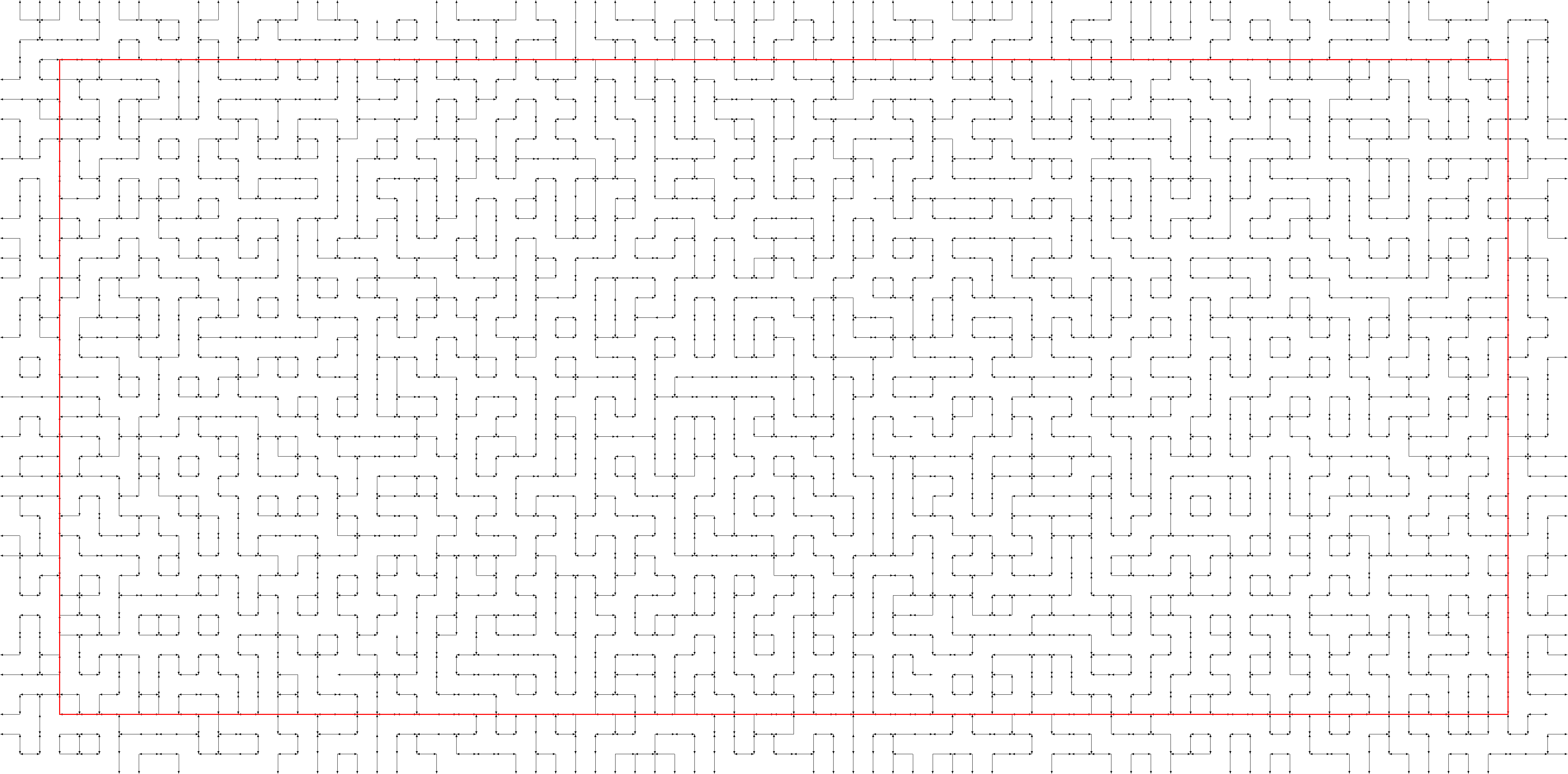}
	\caption{Examples for crossing (left) and non-crossing (right) when $\a = \ff$.}
        \label{fig:sai}
\end{figure}

We start with Theorem \ref{thm:reinf}, i.e., where $\alpha = \infty$. Here, we carried out $N = 10,000$ simulations of the process in a $(80 \times 40)$-rectangle with periodic boundary conditions. In $N_0 = 9,553$ of these simulations we found a horizontal crossing of the central $(72 \times 36)$-rectangle of the nodes that are certainly occupied after $n=4$ steps. Using that the Monte Carlo variance is $\sqrt{0.0447 \cdot 0.9553}$, this implies that with a certainty exceeding $1 - 10^{-300}$, the actual crossing probability is above the threshold of $0.8457$ for 1-dependent percolation The simulations took 22min 17s on a 13th Gen Intel Core i5-1345U. Figure \ref{fig:sai} illustrates examples for crossing and non-crossing realizations.

\begin{figure}[!htpb]
        \centering                                       
\includegraphics[width=0.48\textwidth]{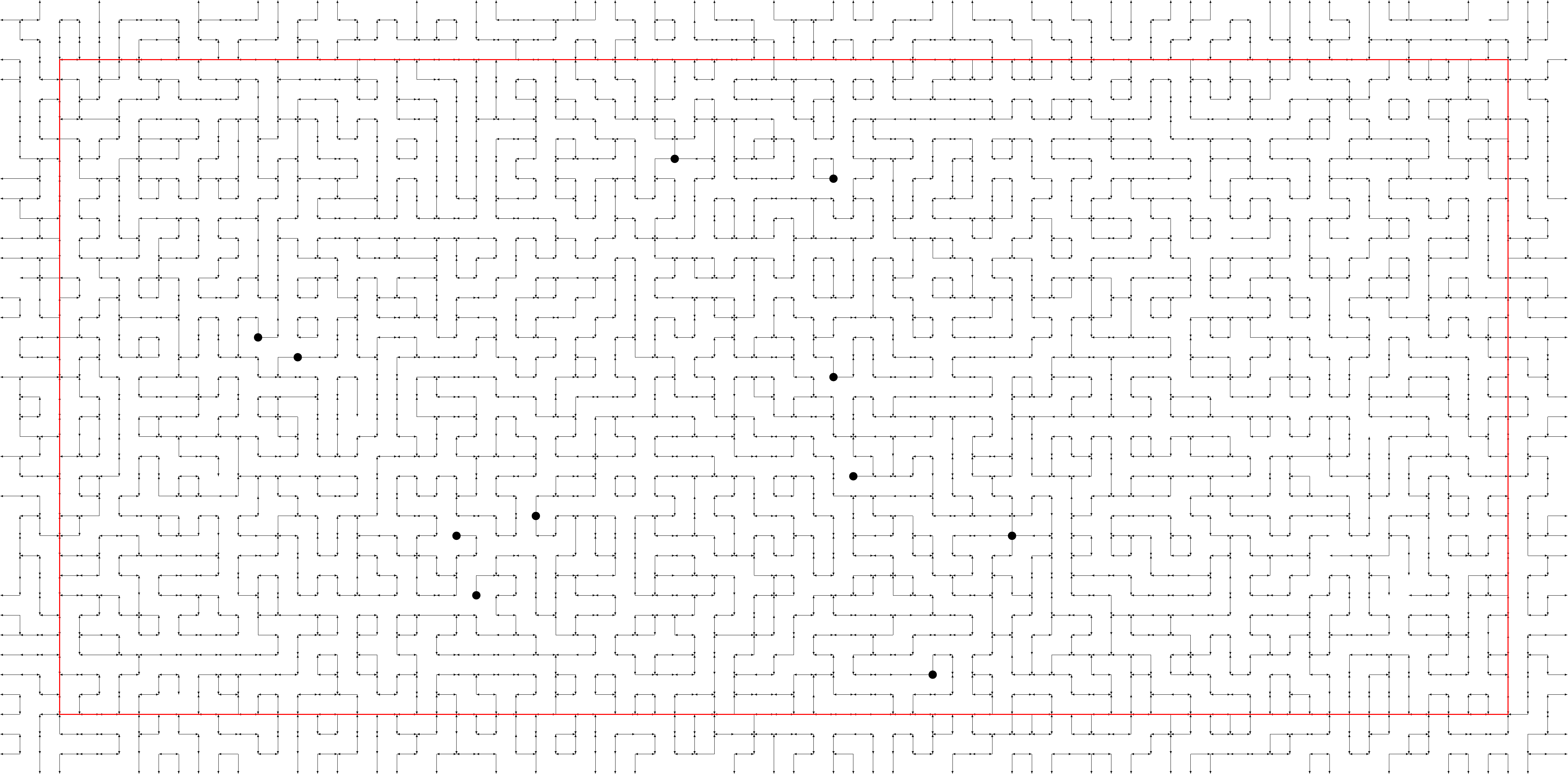}
\includegraphics[width=0.48\textwidth]{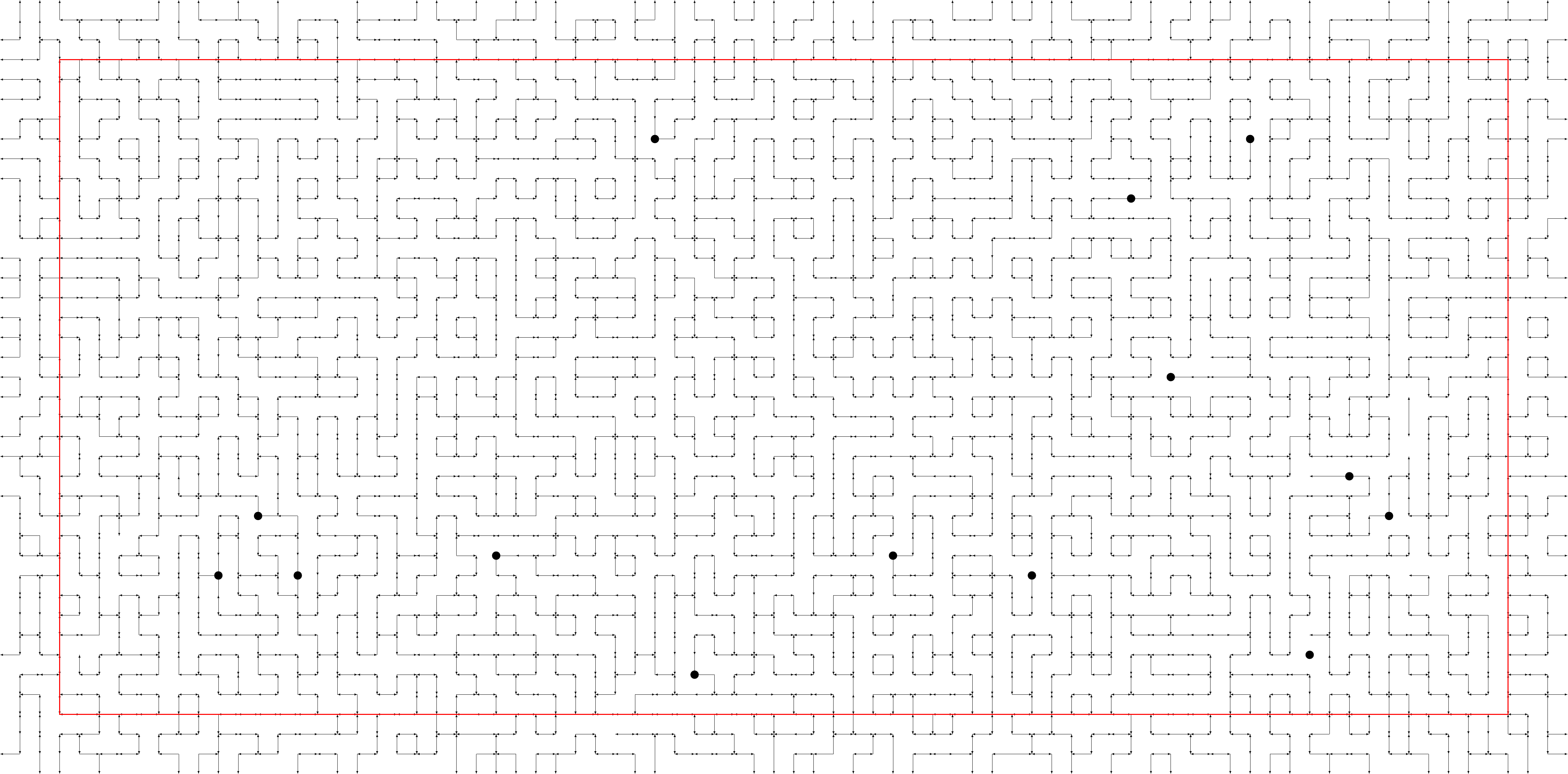}
	\caption{Examples for crossing (left) and non-crossing (right) when $\a = 15$. Corrupted vertices are large black dots.}
        \label{fig:sai2}
\end{figure}

Finally, we discuss Theorem \ref{thm:reinf1}, i.e., where $\alpha < \ff$. Our simulations concern $\a = 15$. The basic setting is the same as above, namely $N = 10,000$ simulations on a $(80 \times 40)$-rectangle with periodic boundary conditions. In $N_0 = 9,512$ of these simulations we found percolation of the central $(72 \times 36)$-rectangle of the nodes that are certainly occupied after $n=4$ steps. Again, with a certainty exceeding $1 - 10^{-300}$, the actual crossing probability is above the threshold of $0.8457$ for 1-dependent percolation. The simulations took 1h 42min 52s on a 13th Gen Intel Core i5-1345U. Note that additional time is needed for sampling from a discrete probability distribution and for dealing with the corruption. Figure \ref{fig:sai2} illustrates examples for crossing and non-crossing realizations.

